\documentclass[letterpaper, 11pt]{article}
\usepackage{amsmath, amsthm, amssymb}    % May not all be necessary
\usepackage{bridges}			% Custom bridges proceedings style
\usepackage{graphicx}			% For including pictures
\usepackage[colorlinks=true, urlcolor=blue, citecolor=black, linkcolor=black]{hyperref}  
\usepackage{subcaption}			% For captioning multi-panel figures
\usepackage{wrapfig}			% wrapping figures with text
\usepackage{comment}
\usepackage[marginpar]{todo}

\newcommand{\R}{\mathbb{R}}
\newcommand{\Nil}{\mathrm{Nil}}
\newcommand{\Sol}{\mathrm{Sol}}
\newcommand{\SL}{\widetilde{\mathrm{SL}_2\mathbb{R}}}

\title{Non-Euclidean Virtual Reality III: Nil}

\author{Remi Coulon\textsuperscript{1}, Elisabetta A. Matsumoto\textsuperscript{2}, Henry Segerman\textsuperscript{3}, Steve Trettel   \textsuperscript{4}
\vspace{10pt}\\
\begin{minipage}[b]{0.5\textwidth}
\centering
\textsuperscript{1}Univ Rennes, CNRS, France; remi.coulon@univ-rennes1.fr\\
\vspace{5pt}
\textsuperscript{3}Oklahoma State University; henry@segerman.org
\end{minipage}
\begin{minipage}[b]{0.5\textwidth}
\centering
\textsuperscript{2}Georgia Institute of Technology; sabetta@gatech.edu\\
\vspace{5pt}
\textsuperscript{4}Stanford University;\\ trettel@stanford.edu
\end{minipage}
}

\date{}					% Suppress any date on submissions

\begin{document}

\maketitle

% Prevent page number 1 from being printed on the first page.
\thispagestyle{empty}

\begin{abstract}

We describe a method of rendering real-time scenes in Nil geometry, and use this to give an expository account of some interesting geometric phenomena.  You can play around with the simulation at \href{www.3-dimensional.space/nil.html}{3-dimensional.space/nil.html}.
\end{abstract}

\begin{figure}[h!tbp]
	\centering
	\includegraphics[width=0.9\textwidth]{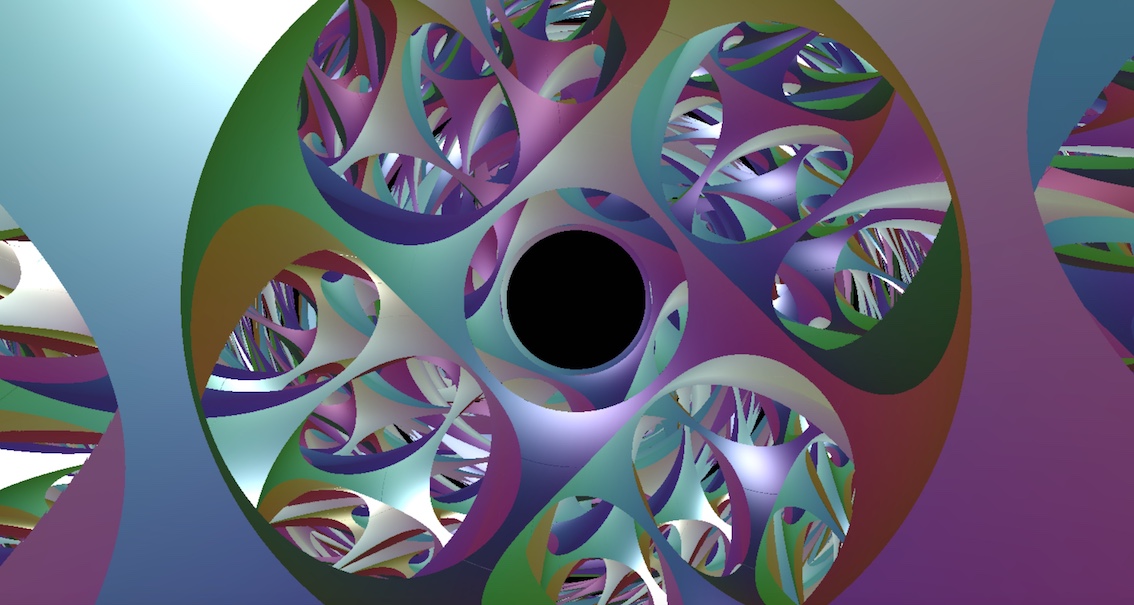}
		\caption{Intrinsic view of a Nil manifold built as a torus bundle over the circle with Dehn twist monodromy.}
	\label{fig:1}
\end{figure}

The Uniformization Theorem of complex analysis tells us that there are three important two-dimensional geometries: hyperbolic, euclidean and spherical.
In contrast, the correct three-dimensional generalization, Thurston's Geometrization Conjecture, singles out eight geometries fundamental to three-dimensional topology.
These come in three flavors: 
\begin{enumerate}
\item constant curvature spaces, $\mathbb{H}^3,\mathbb{E}^3,\mathbb{S}^3$, 
\item products of constant curvature spaces $\mathbb{H}^2\times\R, \mathbb{S}^2\times\R$, and 
\item new geometries built from the three-dimensional Lie groups $\Nil,\Sol,\SL$.
\end{enumerate}
These last three geometries are the most difficult to understand, and a description of their intrinsic geometry is largely missing from the literature.
This paper helps fill this gap by providing an expository account of Nil geometry, 
using the results of a larger project to develop good accurate, real time, intrinsic, and mathematically useful illustrations of homogeneous (pseudo)-riemannian spaces.

\newpage

\section*{Nil Geometry}

%%%%%%%%%%%%%%%%%%%%%%%%%%%%%%%%%%%%%%%%%%%
\begin{wrapfigure}{r}{0.2\textwidth}
\centering
	\includegraphics[width=0.11\textwidth]{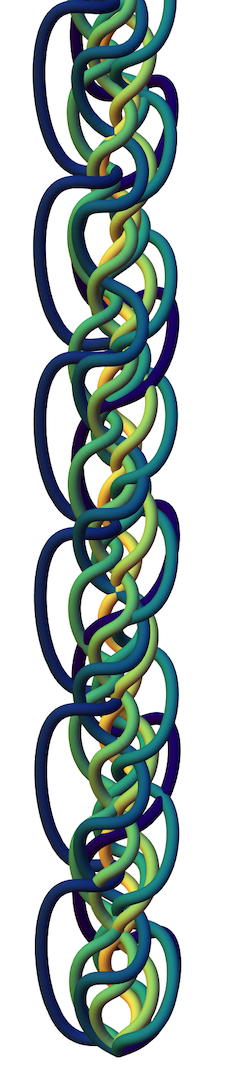}
\caption{Geodesics}
\emph{in the Heisenberg model of Nil coiling around the $z$-axis}
\label{fig:geo in nil}
\end{wrapfigure}

Nil geometry is constructed from the \emph{real Heisenberg group}, which is nilpotent, hence the name.
This is the group of $3\times 3$ matrices of the form
$$\begin{pmatrix}
1&x&z\\0&1&y\\ 0&0&1	
\end{pmatrix}$$
\noindent which we identify with $\R^3$ through the $(x,y,z)$ coordinates.
The center of Nil is the $z$-axis under this identification, which also stands out as a special direction in the geometry.
The metric $ds^2=dx^2+dy^2+(dz-xdy)^2$ is invariant under the action of the Heisenberg group on itself by multiplication, and defines the geometry as $\Nil=(\R^3,ds^2)$.
The metric $ds$ has a non-obvious additional one-parameter family of isometries, which perform a generalized rotation fixing the $z$-axis pointwise.
These rotational isometries are helpful both in calculation, and in explaining some properties of Nil we will encounter along the way.

%\noindent
Restricting the Nil metric to the plane  $\{x=0\}$, we see that the plane is euclidean (but not geodesically embedded). Rotational symmetry then implies that all vertical planes are as well.
%The $xy$ plane in Nil is also copy of $\mathbb{E}^2$ - geodesically embedded in fact - sohorizontal straight lines are Nil geodesics.However, while all horizontal and vertical slices are euclidean, Nil itself is fanother family of euclidean planes in Nil
For each $c \in \R$, the subspace $P_c=\{ z - xy/2 = c\}$ is a totally geodesic subspace isometric to the euclidean plane.
It is invariant under the rotation symmetries on Nil, and its generic geodesics are parabolic arcs.
Nil geodesics outside $P_c$ exhibit more complex behavior and trace out a generalized spiral about the vertical as they travel upwards or downwards, see \autoref{fig:geo in nil}.
%This is quite different from euclidean geometry, as each non-horizontal geodesic in Nil travels only a bounded distance from the $z-$axis.

The geometry of geodesic spheres in Nil depends heavily on their radius.
%,  and displays symmetries coming from the generalized rotations leaving $ds^2_\Nil$ invariant.
Spheres of large radius appear to `dent' inwards on the top and bottom. 
This signals a sort of difficulty in traveling vertically, and indeed it is often \emph{quicker} to reach a point $(0,0,c)$ by following a wide helix rather than traveling along the $z-$axis.

 \begin{figure}[h!tbp]
\centering
 	\includegraphics[width=0.9\textwidth]{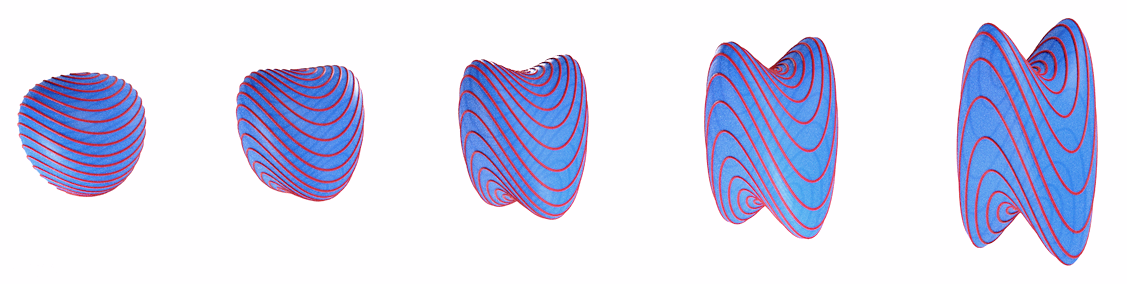}
    \caption{Balls of radius $1$, $2$, $3$, $4$, and $5$ in the Heisenberg model of Nil.
The balls have been rescaled to have approximately the same size.
The curves drawn on the balls are invariant under the rotational symmetries of Nil.}
\label{fig:2}
\end{figure}

\vspace{-0.5cm}

\section*{Seeing Inside Nil Geometry}
 
What it would actually look like if you woke up one day inside of Nil geometry?
It helps to think briefly about the mechanics of sight.
The image formed upon our retinas is created from light rays bouncing off objects in the world and into our eyes.
Along the path from object to eye, light follows a \emph{geodesic} -- in everyday life, a straight line\footnote{  
On cosmic scales this is not true even in our world, as evidenced by telescopic images of \emph{gravitationally lensed} distant galaxies.}.
Thus, to model sight in Nil, we place an imaginary screen in space -- our simulated retina -- and trace outwards along Nil geodesics from this screen to objects in our simulated scene.

\begin{figure}[h!tbp]
\centering
\begin{minipage}[b]{0.22\textwidth} 
	\includegraphics[width=\textwidth]{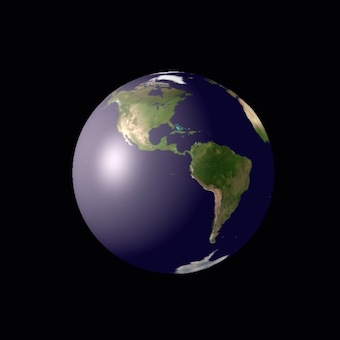}%
        	\subcaption{} % Add subcaption text if desired, or use \subcaption* to suppress (a), (b), etc. labels
        	\label{fig:rings12a}
\end{minipage}
~ %add desired spacing between images, e. g., ~, \quad, \qquad, \hfill etc.	
\begin{minipage}[b]{0.22\textwidth} 
	\includegraphics[width=\textwidth]{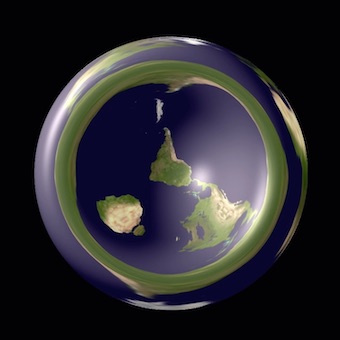}%
        	\subcaption{} % Add subcaption text if desired, or use \subcaption* to suppress (a), (b), etc. labels
        	\label{fig:rings12b}
\end{minipage}
~ %add desired spacing between images, e. g., ~, \quad, \qquad, \hfill etc.	
\begin{minipage}[b]{0.22\textwidth} 
	\includegraphics[width=\textwidth]{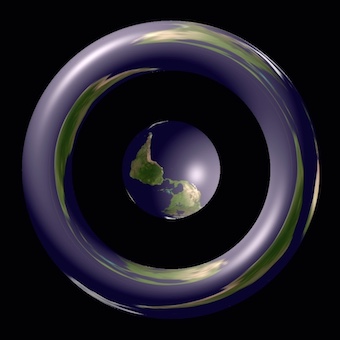}%
        	\subcaption{} % Add subcaption text if desired, or use \subcaption* to suppress (a), (b), etc. labels
        	\label{fig:rings12c}
\end{minipage}
~ %add desired spacing between images, e. g., ~, \quad, \qquad, \hfill etc.	
\begin{minipage}[b]{0.22\textwidth} 
	\includegraphics[width=\textwidth]{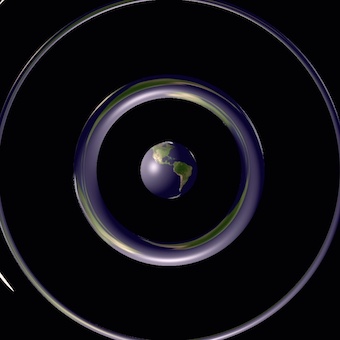}%
        	\subcaption{} % Add subcaption text if desired, or use \subcaption* to suppress (a), (b), etc. labels
        	\label{fig:rings12d}
\end{minipage}
\caption{Flying away from a sphere in Nil, along the $z$-axis (a) before the first conjugate point,
(b,c) past the first conjugate point, before the second conjugate point, \& (d) past the second conjugate point. }
\label{fig:rings1}
\end{figure}

To highlight some geometric features in Nil, we focus on a simple scene: an otherwise empty universe containing a single sphere (textured as the earth).
Figure \ref{fig:rings1} shows successive views from a spaceship flying away from the earth along the $z$-axis.
They tell a rather surprising, and perhaps confusing story.
While at first the earth appears normal, it quickly becomes distorted and then appears to slough off concentric \emph{rings} as it recedes further into the distance. This is best understood with a video\footnote{\url{https://vimeo.com/382647580}} of the effect.

\begin{figure}[h!tbp]
\centering
	\includegraphics[width=0.85\textwidth]{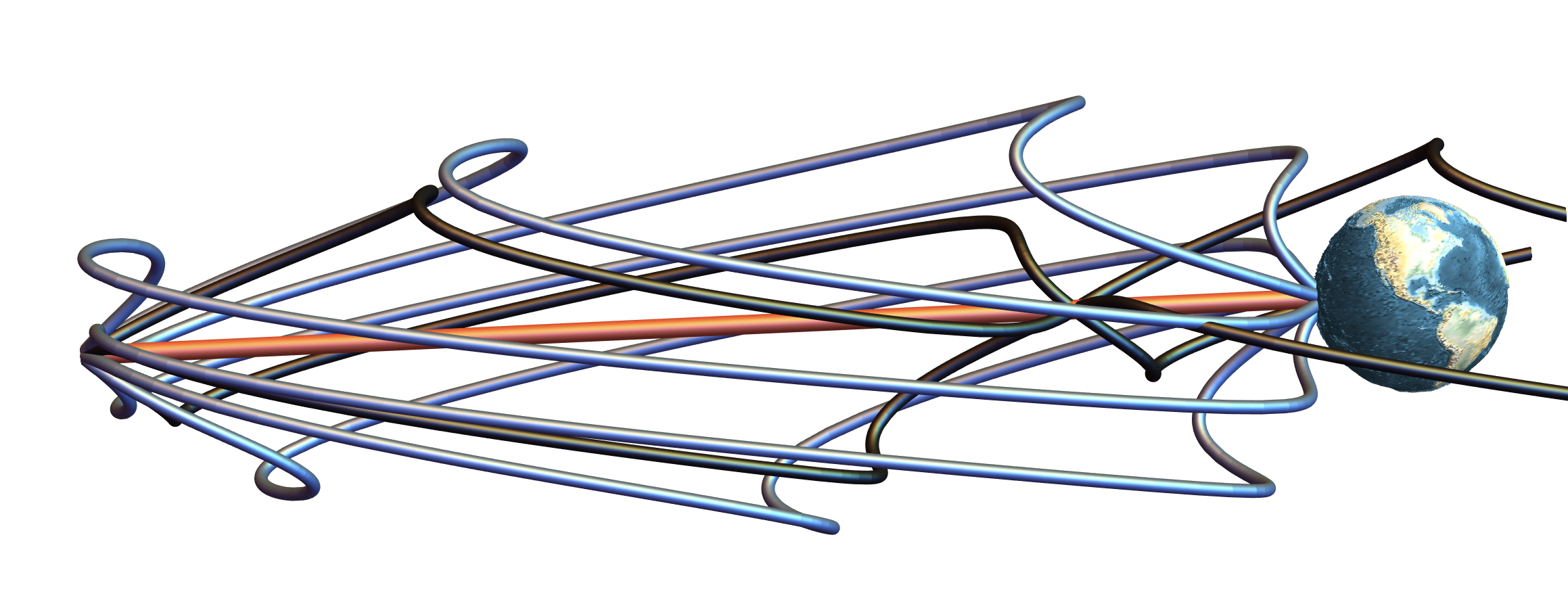}%
        	\label{fig:2a}
\caption{The cause of the ring-shaped mirages.
The inner geodesic (red) and ring of outer geodesics (blue) both reach the same point on earth.  In between, the black geodesics escape into space.}
\label{fig:extrinsic1}
\end{figure}

\noindent
%Recall that points on the screen are colored by \emph{following geodesics}.
A pair of points $p,q$ are \emph{conjugate} if $p$ is joined to $q$ by multiple geodesics. 
These appear along the $z$-axis in Nil, for pairs of points that are sufficiently far apart. In \autoref{fig:extrinsic1}, we see multiple copies of the same point on the earth -- this happens when the camera and that point are a conjugate pair.
%If you are located at $p$, and $q$ lies on the surface of the earth, then you will see multiple copies of $q$, see .
As Nil has rotational symmetries fixing the $z$-axis, there is actually a continuum of such copies forming an \emph{entire circle of images}.
This \emph{geometric lensing} in Nil is analogous to the gravitational lensing seen around massive objects in the cosmos.

\vspace{-0.5cm}

\noindent

\begin{figure}[h!tbp]
\centering
\begin{minipage}[b]{0.65\textwidth} 
	\includegraphics[width=0.9\textwidth]{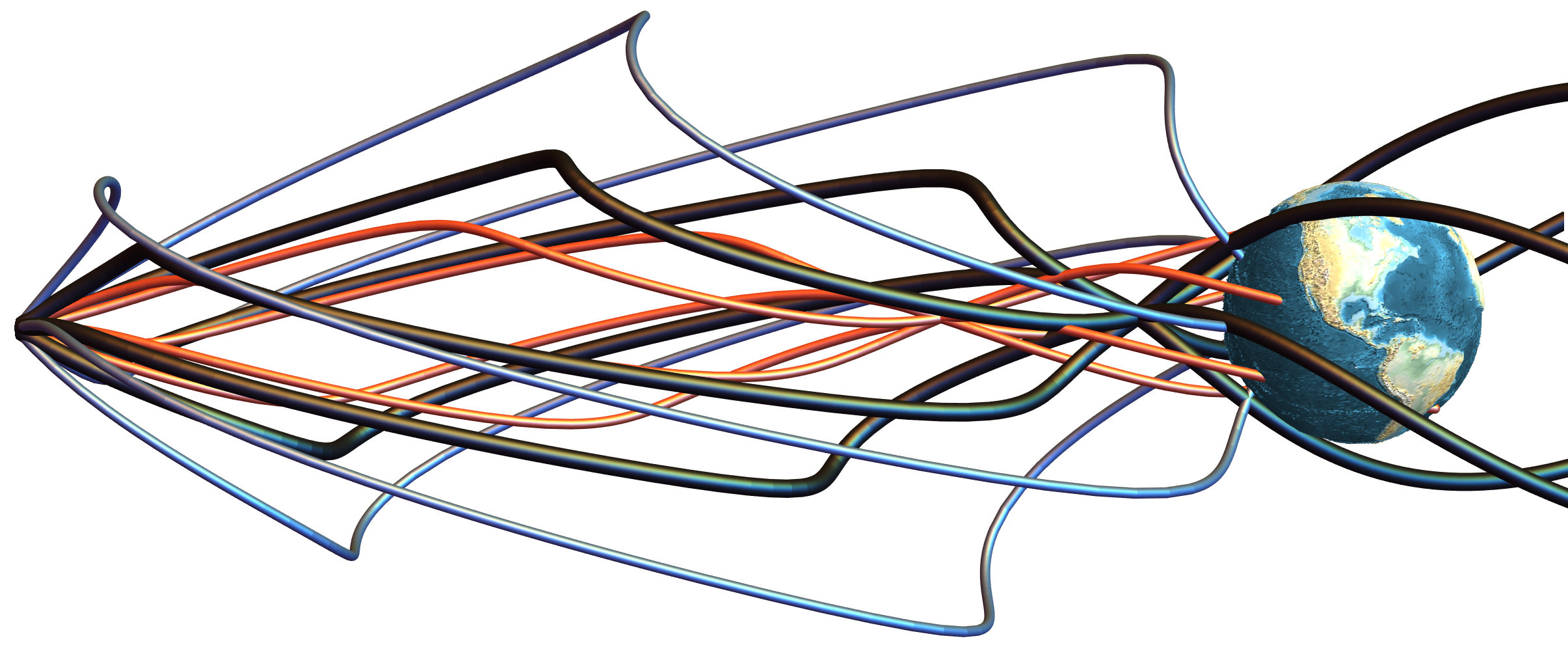}%
\end{minipage}\quad
\begin{minipage}[b]{0.3\textwidth} 
	\raisebox{1cm}{\includegraphics[width=\textwidth]{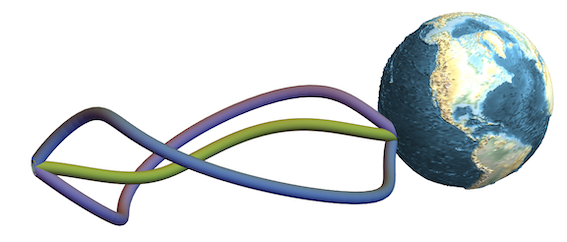}}%
\end{minipage}

\vspace{-1mm}
\caption{
(a) An interior horizon (black geodesics) of tangent geodesics has all nearby geodesics intersecting the earth.  (b) Point on earth with three images: one near the visual center (green), and two farther out, in opposite directions (blue, purple)}
\label{fig:extrinsic2}
\end{figure}

\begin{figure}[h!tbp]
\centering
\begin{minipage}[b]{0.22\textwidth} 
	\includegraphics[width=\textwidth]{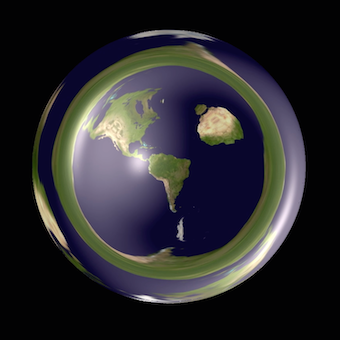}
        	%\subcaption{} % Add subcaption text if desired, or use \subcaption* to suppress (a), (b), etc. labels

\end{minipage}\quad
 %add desired spacing between images, e. g., ~, \quad, \qquad, \hfill etc.	
\begin{minipage}[b]{0.22\textwidth} 
	\includegraphics[width=\textwidth]{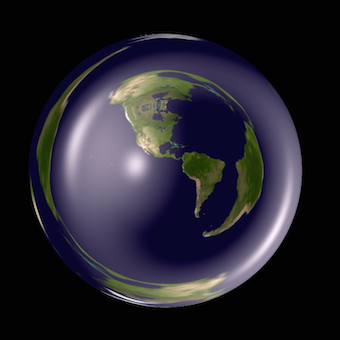}
        	%\subcaption{} % Add subcaption text if desired, or use \subcaption* to suppress (a), (b), etc. labels
\end{minipage}\quad
 %add desired spacing between images, e. g., ~, \quad, \qquad, \hfill etc.	
\begin{minipage}[b]{0.22\textwidth}
	\includegraphics[width=\textwidth]{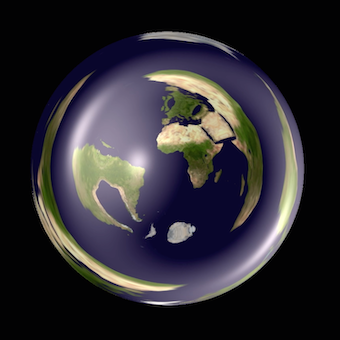}
        	%\subcaption{} % Add subcaption text if desired, or use \subcaption* to suppress (a), (b), etc. labels
\end{minipage}\quad
\begin{minipage}[b]{0.22\textwidth} 
	\includegraphics[width=\textwidth]{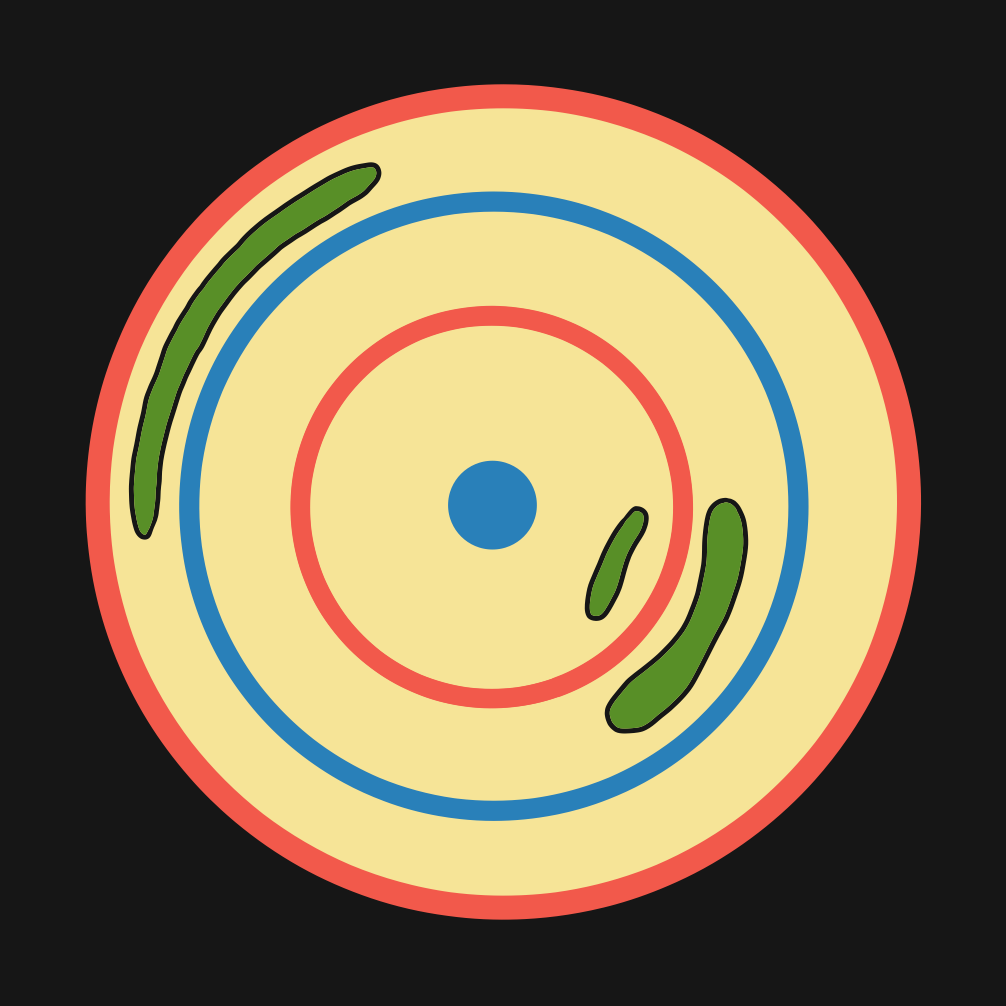}
        	%\subcaption{} % Add subcaption text if desired, or use \subcaption* to suppress (a), (b), etc. labels
\end{minipage}

\caption{\href{https://vimeo.com/382648346}{(Video Link)} Three views of Earth past the first conjugate point.  (Right) A Schematic showing }
\emph{important  features: exact center \& mirage ring (blue), inner \& outer horizons (red), generic region (green).}
\label{fig:cong2}
\end{figure}

Another important geometric feature is the \emph{horizon}, realized by geodesics tangent to the earth.
This has one image in \ref{fig:rings1}a, three images in Figure \ref{fig:rings1}c, and generally gains two more each time a new ring appears.  But
before the imaged earth even splits off a first ring, an \emph{interior horizon} develops. 
This is simulated by the black geodesics in Figure \ref{fig:extrinsic2}a and visible in all of Figure \ref{fig:cong2}: nearby geodesics both inside and outside this tangency envelope intersect the earth, causing an image of the horizon surrounded by earth on both sides.

Figure \ref{fig:extrinsic2}b analyzes the image of points on the earth which are not directly on the axis of symmetry.
Before the second conjugate point, generic points on the earth's surface have three images: one from a near-central geodesic making one full twist, and more coming from geodesics each completing half a turn.
The innermost image is separated from the second by an image of the horizon, and the second from the third by the circle's worth of additional images of the center.

\vspace{-0.25cm}

\noindent 
Figure \ref{fig:cong2} shows all of these effects in a sequence of images, depicting the earth as it rotates.
The first image is centered South America, which then is also imaged as a large green ring. 
Next this unstable ring breaks into two new, opposite images, one `inside' and one `outside' the original ring.
The third image makes the interior horizon visible as we see two copies of South America beginning to disappear into it.

\noindent 
The instability of the ring-shaped South America in Figure \ref{fig:cong2} applies to the entire ring-shaped mirage of the earth as well -- which causes real difficulty for creatures evolved to perceive depth through parallax such as ourselves.
Figure~\ref{fig:2} shows a sequence of views of earth with only a small horizontal change in position -- such as the separation between our eyes --  from one to the next.
Near conjugate points, the images our brain receives from each eye may have very little in common!

\begin{figure}[h!tbp]
\centering
\begin{minipage}[b]{0.18\textwidth} 
	\includegraphics[width=\textwidth]{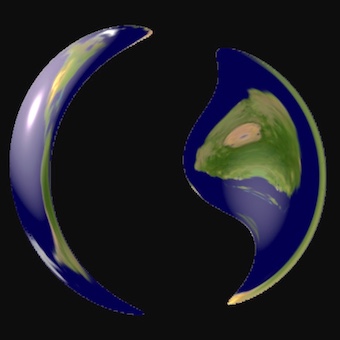}
\end{minipage}
~ %add desired spacing between images, e. g., ~, \quad, \qquad, \hfill etc.	
\begin{minipage}[b]{0.18\textwidth} 
	\includegraphics[width=\textwidth]{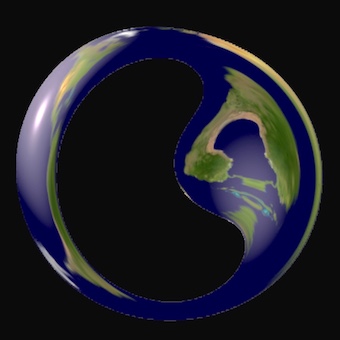}
\end{minipage}
~ %add desired spacing between images, e. g., ~, \quad, \qquad, \hfill etc.	
\begin{minipage}[b]{0.18\textwidth} 
	\includegraphics[width=\textwidth]{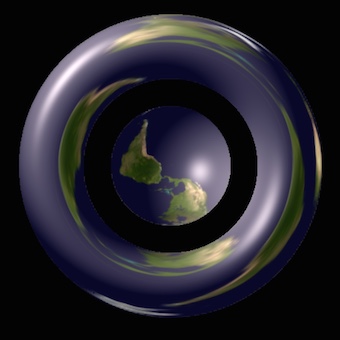}
\end{minipage}
~ %add desired spacing between images, e. g., ~, \quad, \qquad, \hfill etc.	
\begin{minipage}[b]{0.18\textwidth} 
	\includegraphics[width=\textwidth]{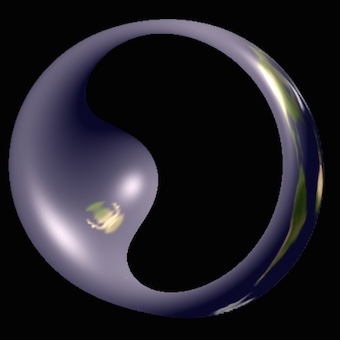}
\end{minipage}
~
\begin{minipage}[b]{0.18\textwidth} 
	\includegraphics[width=\textwidth]{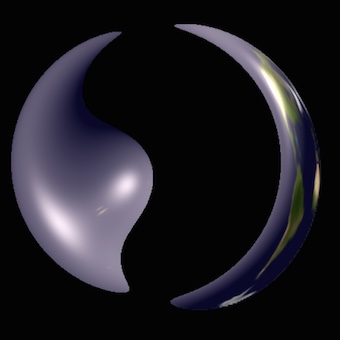}

\end{minipage}
\caption{Parallax in $\mathsf{Nil}$ makes stereoscopic vision difficult.  Images are rendered at a distance of $\sim 10$ units from the earth, shifting horizontally by approximately $\sim 0.1$ per frame. \href{https://vimeo.com/385915362}{(Video Link)}}
\label{fig:2}
\end{figure}

\noindent
Another surprising fact: 
in the sequence of images in Figure \ref{fig:rings1}, the earth appears bigger and bigger even though we are moving further away from it.
Indeed the angle between the outside geodesics reaching the earth increases as we fly away from the earth.
No matter how far away the earth gets, the angular size of its central image remains unchanged -- a consequence of the fact that geodesics in Nil only travel a bounded distance from the $z$ axis.

\begin{figure}[h!tbp]
\centering
\includegraphics[width=0.9\textwidth]{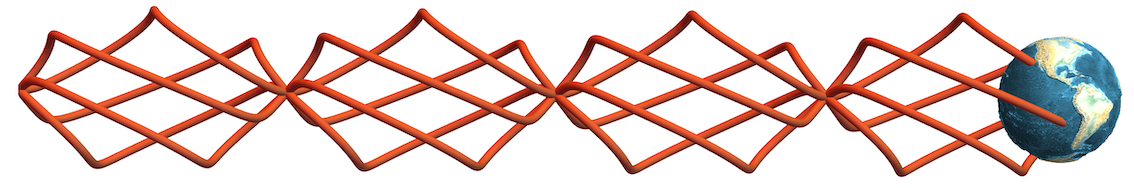}
\vspace{-1mm}
\caption{
Even far away, Earth takes up a large portion of the visual field as geodesics stay near the $z$-axis.}
\label{fig:extrinsic3}
\end{figure}

\begin{figure}[h!tbp]
\centering
\begin{minipage}[b]{0.2\textwidth} 
	\includegraphics[width=\textwidth]{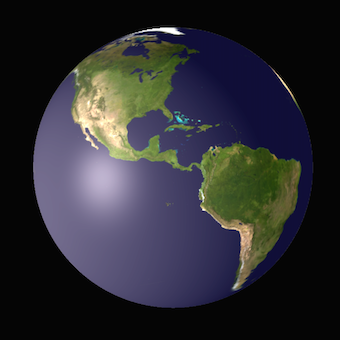}
\end{minipage}\quad
\begin{minipage}[b]{0.2\textwidth} 
	\includegraphics[width=\textwidth]{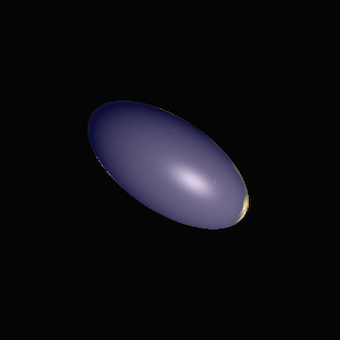}
\end{minipage}\quad
\begin{minipage}[b]{0.2\textwidth} 
	\includegraphics[width=\textwidth]{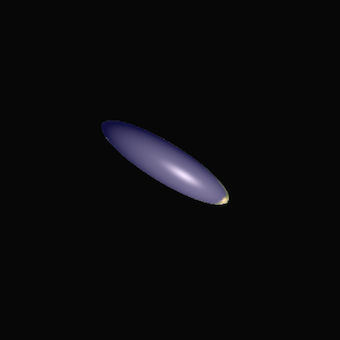}
\end{minipage}\quad
\begin{minipage}[b]{0.2\textwidth} 
	\includegraphics[width=\textwidth]{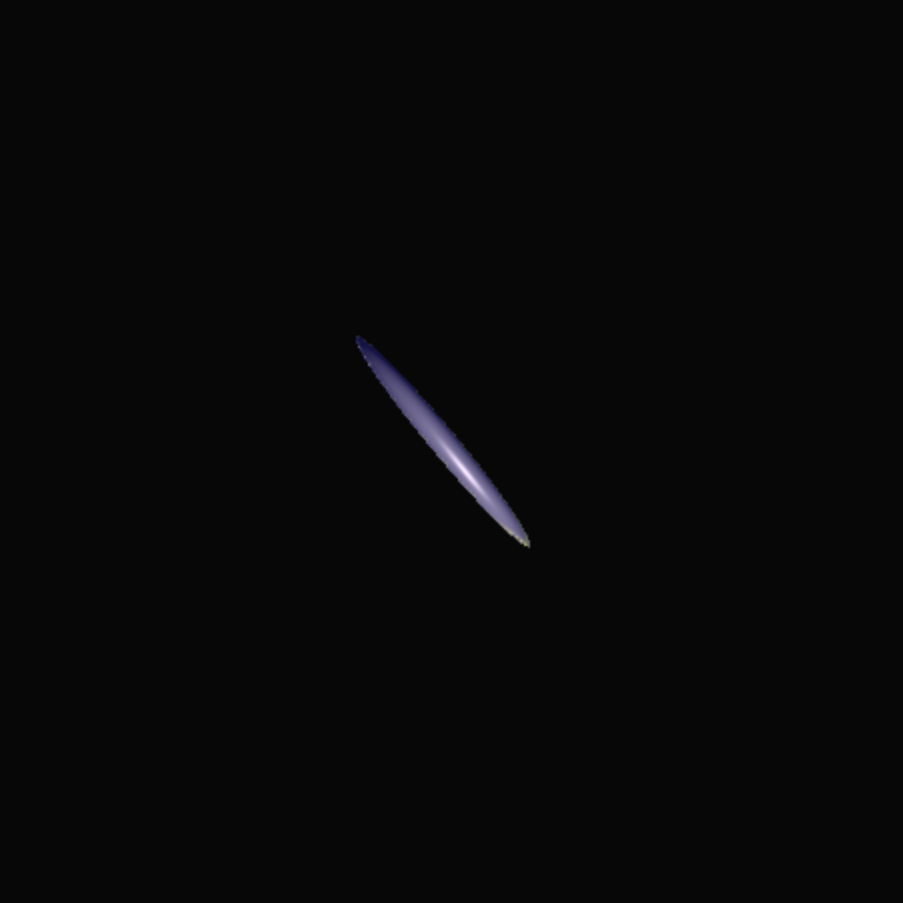}
\end{minipage}
\caption{Flying away from a sphere in Nil, along the $x$-axis.}
\label{fig: away along x}
\end{figure}

\noindent
We have focused on imaging a sphere from points along the $z$-axis, as this best shows the behavior of conjugate points in Nil.
To close this section we offer one more perspective: flying away from a sphere while staying in the euclidean subspace $\{ z - xy/2 = 0\}$, see \autoref{fig: away along x}.
There are no conjugate points to our location in these directions, and the earth shrinks as we recede from it.
Note however that distant earths appear greatly elongated and rotated.
This behavior exemplifies the anisotropy of Nil. 

\vspace{-0.25cm}

\section*{Inside Compact Nilmanifolds}

\begin{figure}[h!tbp]
\centering
\begin{minipage}[b]{0.22\textwidth} 
	\includegraphics[width=\textwidth]{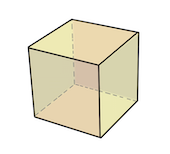}
        	%\subcaption{} % Add subcaption text if desired, or use \subcaption* to suppress (a), (b), etc. labels

\end{minipage}
~ %add desired spacing between images, e. g., ~, \quad, \qquad, \hfill etc.	
\begin{minipage}[b]{0.24\textwidth} 
	\includegraphics[width=\textwidth]{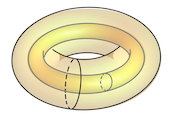}
        	%\subcaption{} % Add subcaption text if desired, or use \subcaption* to suppress (a), (b), etc. labels
\end{minipage}
~ %add desired spacing between images, e. g., ~, \quad, \qquad, \hfill etc.	
\begin{minipage}[b]{0.4\textwidth} 
	\includegraphics[width=\textwidth]{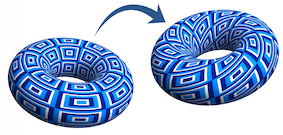}
        	%\subcaption{} % Add subcaption text if desired, or use \subcaption* to suppress (a), (b), etc. labels
\end{minipage}
\caption{Identifying vertical faces of the cube in opposite pairs gives a thickened torus, and   gluing these boundary tori via a Dehn Twist (right) produces a closed Nil manifold.}
\label{fig:nil manifold}
\end{figure}

Closed manifolds locally modeled on Nil geometry are some of the fundamental building blocks for three-dimensional topology.
 Such manifolds are covered by \emph{Dehn-twist torus bundles}, which have a nice topological construction. See \autoref{fig:nil manifold}.
 Starting with a cube, identify opposing pairs of vertical sides so that each cross section becomes a torus.
 The resulting space is a \emph{thickened torus}. To create a closed three-manifold we need only glue the inner surface to the outer surface.
 This manifold has Nil geometry when we glue with a 
\emph{Dehn twist}: the result of cutting a torus along a closed curve, twisting a full $2\pi n$, and re-gluing, as seen in the final illustration in \autoref{fig:nil manifold}.

\begin{figure}[h!tbp]
\centering
\begin{minipage}[b]{0.32\textwidth} 
	\includegraphics[width=\textwidth]{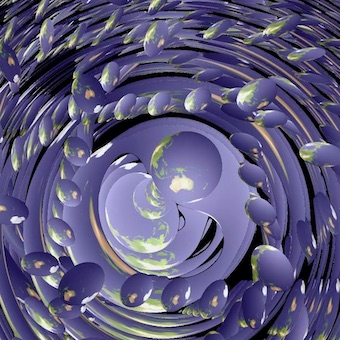}%
        	\subcaption{} % Add subcaption text if desired, or use \subcaption* to suppress (a), (b), etc. labels
        	\label{fig:2a}
\end{minipage}
~ %add desired spacing between images, e. g., ~, \quad, \qquad, \hfill etc.	
\begin{minipage}[b]{0.32\textwidth} 
	\includegraphics[width=\textwidth]{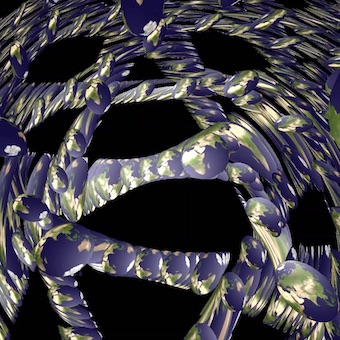}%
        	\subcaption{} % Add subcaption text if desired, or use \subcaption* to suppress (a), (b), etc. labels
        	\label{fig:2b}
\end{minipage}
~ %add desired spacing between images, e. g., ~, \quad, \qquad, \hfill etc.	
\begin{minipage}[b]{0.32\textwidth} 
	\includegraphics[width=\textwidth]{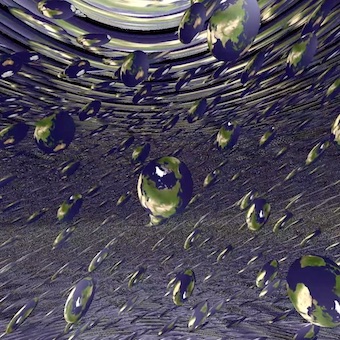}%
        	\subcaption{} % Add subcaption text if desired, or use \subcaption* to suppress (a), (b), etc. labels
        	\label{fig:2c}
\end{minipage}
\caption{\href{https://vimeo.com/385893478}{(Video Link)} An earth in the Nilmanifold formed as a Dehn twist torus bundle. \\
(a) Facing along the z-axis,  (b) diagonal view,  (c) looking in the direction of the xy plane}
\label{fig:EarthLattice}
\end{figure}

\noindent 
The images shown in Figure \ref{fig:EarthLattice} are the result of placing a single earth in the torus bundle manifold described above.
The numerous images of earth have two sources: light getting wrapped around
the manifold's nontrivial topology, and the \emph{geometric lensing} effects we saw above.
The rings in Figures \ref{fig:EarthLattice}a,c are geometrically lensed images of some far-away earth lined up with us along the $z$-axis, and the multitude of foreground earths visible in Figures \ref{fig:EarthLattice}b,c are from light wrapping around nontrivial loops in the torus bundle.

\begin{figure}[h!tbp]
\centering
\begin{minipage}[b]{0.32\textwidth} 
	\includegraphics[width=\textwidth]{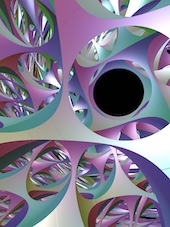}%
        	\subcaption{} % Add subcaption text if desired, or use \subcaption* to suppress (a), (b), etc. labels
        	\label{fig:2a}
\end{minipage}%\quad
~ %add desired spacing between images, e. g., ~, \quad, \qquad, \hfill etc.	
\begin{minipage}[b]{0.32\textwidth} 
	\includegraphics[width=\textwidth]{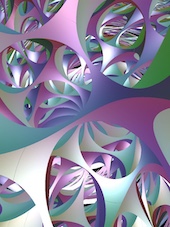}%
        	\subcaption{} % Add subcaption text if desired, or use \subcaption* to suppress (a), (b), etc. labels
        	\label{fig:2b}
\end{minipage}%\quad
~ %add desired spacing between images, e. g., ~, \quad, \qquad, \hfill etc.	
\begin{minipage}[b]{0.32\textwidth} 
	\includegraphics[width=\textwidth]{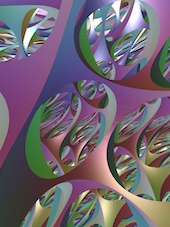}%
        	\subcaption{} % Add subcaption text if desired, or use \subcaption* to suppress (a), (b), etc. labels
        	\label{fig:2c}
\end{minipage}
\caption{\href{https://vimeo.com/385890730}{(Video Link)} Edges of a fundamental domain for the same Nil manifold as above.  \\
(a) Facing along the z-axis,  (b) diagonal view,  (c) looking in the direction of the xy plane}
\label{fig:tiling}
\end{figure}

For another perspective, we recall that this manifold was built from a cubical fundamental domain and render its edges. Multiple images of these produce a tiling of Nil shown in Figure \ref{fig:tiling}.
Being more `full of stuff' actually makes this scene somewhat easier to navigate, as foreground multiple-images coming from the nontrivial topology obscure many of the more confusing effects of geometric lensing in the background.

%%%%%%%%%%%%%%%%%%%%%%%%%%%%%%%%%%%%%%%
\section*{Implementation Details}
These intrinsic simulations of Nil were created by adapting the computer-graphics techinque of \emph{ray-marching} to non-euclidean homogeneous spaces.
This differs from previous work of the authors \cite{NEVR1,NEVR2} which computes intrinsic views by pulling objects back to the tangent space via partial inverses of the riemannian exponential map.
This method is not well-adapted to the geometry of $\Nil$ because the exponential map is neither invertible nor has an easily computable cut-locus.
Instead of attempting to pull things back along geodesics, from each pixel on the screen we flow out along a geodesic out into the space. Upon intersecting an object we compute its color by then following geodesics to the light sources then apply the Phong reflection model~\cite{phong}.
The precise implementation details are the subject of a forthcoming paper, and the code (currently a work in progress) is available on GitHub \cite{github}.
Other simulations of Nil geometry include the work of Berger \cite{berger}, and ZenoRogue (who is producing games in the Thurston geometries!) \cite{Zeno2,ZenoR}.

\vspace{2mm}
\noindent
\emph{\textbf{\large Computing in Nil:\;}}
The Heisenberg model of the Nil geometry discussed above is very easy to state but has a major drawback: the rotational symmetries of Nil do not correspond to the euclidean rotation around the $z$-axis.
In practice we work with a rotation invariant set of coordinates. The map from the Heisenberg model to the rotation invariant one is given by $(x,y,z) \to (x,y,z - xy/2)$.
%The metric in these new coordinates is GIVE METRIC.
We store our location in Nil by a pair consisting of a point $p\in\Nil$ and an orthonormal frame $F=(e_1,e_2,e_3)$ based at $p$.  
The frame determines our orientation: the image rendered on the screen is rendered so that we are looking straight down the $-e_3$ axis, with $e_2$ to the top and $e_1$ to the right.
Since the underlying space of Nil is the Heisenberg group, the point $p$ also uniquely specifies an \emph{isometry $L_p$ taking the origin $o = (0,0,0)$ to $p$}. To simplify computation we store not $F$ but its translate to the origin, $dL_p^{-1}(F)$; a matrix in $\mathsf{SO}(3)$.

\vspace{2mm}
\noindent
\emph{\textbf{\large Producing a Fixed Image:\;}}
To compute the intrinsic view of a scene $S$ from a position $(p,F)$, we compute from our position the collection of tangent vectors representing pixels on the screen. To continue we require a means of computing the geodesic flow in Nil.
Given a unit length vector $v=(a,0,c)\in T_0\Nil$, the geodesic equation $\nabla_{\dot{\gamma}}\dot{\gamma}=0$ for $\gamma(0)=0$, $\dot{\gamma}(0)=v$ can be solved explicitly. This determines all geodesics through the translational and rotational symmetries of Nil.
\vspace{-1mm}
\begin{equation*}
    \gamma(t) =\left(
    \frac{2a}c\sin (c t)\cos(c t),\
    \frac{2a}c\sin ^2(c t),\
    c t+\frac{a^2}{2c^2}\left[c t-\sin(c t)\right]
    \right).
\end{equation*}

\noindent
To calculate how far to flow and to determine when we have reached an object, we need a \emph{signed distance function} measuring the distance in Nil from each point to the nearest object in the scene $S$. For example, the  distance function from a point $p = (x_1, x_2, x_3)$ to a sphere of radius $r$ centered at the origin is $d(p,o)-r$, where $d$ is the distance in Nil.
Unfortunately this is difficult to compute, so we approximate as follows.
We estimate a function $F \colon \Nil \to \R_+$ which bilipschitz approximates $d(p,o)$, and use this as a good-enough-estimate when far away.
For two nearby points, constructing the minimizing geodesic joining them is numerically feasible. So once $F$ passes a certain threshhold, we switch to numerically calculating the actual distance.

\vspace{2mm}
\noindent
\emph{\textbf{\large Moving in Nil \;}}
To change one's facing in Nil, we act by rotation on the frame which stores orientation, without changing position.
%This changes the initial directions of the unit tangent vectors for the geodesic flow without changing their basepoint.
To change position, we take as input (from the keyboard, headset, etc.) a displacement vector $v\in\R^3$, and interpret it as a tangent vector at our current position $p\in\Nil$.
We follow the geodesic flow by $v$ from $p$; update our position to its end point $q$ and
update our facing by the parallel transport associated to this geodesic.

\vspace{2mm}
\noindent
\emph{\textbf{\large Quotient Manifolds:\;}}
The closed Nil-manifold $M$ depicted in this paper is the mapping torus of the Dehn twist corresponding to the matrix $A=\left(\begin{smallmatrix}1&1\\0&1\end{smallmatrix}\right)$.
To ray-march in $M$ we need both a fundamental domain and an algorithm for `teleporting' back into the fundamental domain to continue the geodesic flow.
A straightforward calculation shows we may take the unit cube in Nil as a fundamental domain, with side pairings $L_{(1,0,0)}$, $L_{(0,1,0)}$ and $L_{(0,0,1)}$.
To return an escaped point back to the fundamental domain, we follow a simple-minded algorithm: if the point in question has $x\not\in[0,1]$, we iteratively apply the $x$ shear-translation generator until $x\in[0,1]$.
Then, in either order (as they commute) we iteratively apply the $y$ and $z$ translations (which do not affect the $x$ coordinate) until the point returns to the unit cube.

%%%%%%%%%%%%%%%%%%%%%%%%%%%%%%%%%%%%%%%

%%%%%%%%%%%%%%%%%%%%%%%%%%%%%%%%%%%%%%%
\section*{Summary and Future Work}

This project has produced a real-time, intrinsic and geometrically correct rendering engine for Nil geometry and its compact quotients, that can take movement input from either a keyboard or headset, and can render images from two viewpoints simultaneously (for stereoscopic vision).
However, this is still very much a work in progress.
Currently we aim to calculate the signed distance functions for more basic geometric objects in Nil (including, for instance, the distance to a vertical or horizontal plane) and to improve our lighting calculations (which right now, only account for the light that hits an object following the \emph{shortest} geodesic connecting the light to that object, not \emph{all} geodesics).

%%%%%%%%%%%%%%%%%%%%%%%%%%%%%%%%%%%%%%%

%%%%%%%%%%%%%%%%%%%%%%%%%%%%%%%%%%%%%%%
\section*{Acknowledgements}
This material is based upon work supported by the National Science Fundation under Grant No. DMS-1439786 while the authors were in residence at the ICERM in Providence, RI, during the semester program Illustrating Mathematics.
The first author also acknowledges support from the Agence Nationale de la Recherche under Grant \emph{Dagger} ANR-16-CE40-0006-01 as well as the \emph{Centre Henri Lebesgue} ANR-11-LABX-0020-01. The second author is grateful to support from the National Science Foundation DMR-1847172. The third author was supported in part by National Science Foundation grant DMS-170823.
We are especially thankful to Brian Day for helping get things synced up with the virtual reality headsets, and are excited to incorporate some of his work on non-euclidean physics in the future.
This project is indebted to a long history of previous work.  It is a direct descendant of the hyperbolic ray-march program created by Roice Nelson, Henry Segerman, and Michael Woodard~\cite{RSC}, which itself was inspired by previous work in $\mathbb{H}^3$ by Andrea Hawksley, Vi Hart, Elisabetta A. Matsumoto and Henry Segerman~\cite{NEVR1, NEVR2}, all of which aim to expand upon the excellent work of Jeff Weeks in \emph{Curved Spaces}~\cite{Weeks}.

%%%%%%%%%%%%%%%%%%%%%%%%%%%%%%%%%%%%%%%

%%%%%%%%%%%%%%%%%%%%%%%%%%%%%%%%%%%%%%%
% References %
    
{\setlength{\baselineskip}{13pt} % tighten line spacing for bibliography
\raggedright				% no right justification for References

} % end setlength, raggedright

\todos
\end{document}